\documentclass[preprint,3p,sort]{elsarticle}

\usepackage{amsfonts}
\usepackage{amsthm}
\usepackage{graphics}
\usepackage{epstopdf}
\usepackage{algorithmic}
\usepackage{ams math}
\usepackage{amssymb}
\usepackage{hyperref}
\usepackage[T1]{fontenc}
\usepackage{soul,color}
\usepackage{cases}
 \usepackage{lineno}
\newtheorem{theorem}{Theorem}[section]
\newtheorem{proposition}[theorem]{Proposition}
\newtheorem{lemma}[theorem]{Lemma}
\newtheorem{corollary}[theorem]{Corollary}
\newtheorem{definition}{Definition}[section]
\newtheorem{remark}[theorem]{Remark}
\theoremstyle{example}
\newtheorem{example}[theorem]{Example}
\numberwithin{equation}{section}
\begin{document}


\begin{frontmatter}

%




\title{Unconditionally stable fully discrete high order approximations of  \\ nonlinear initial boundary value problems}

\author[sweden,southafrica]{Jan Nordstr\"{o}m}
\address[sweden]{Department of Mathematics, 
Link\"{o}ping University, SE-581 83 Link\"{o}ping, Sweden}
\address[southafrica]{Department of Mathematics and Applied Mathematics, University of Johannesburg, P.O. Box 524, 
Auckland Park 2006, Johannesburg, South Africa
}

\begin{abstract}
Summation-by-parts (SBP) operators combined with simultaneous-approximation-terms (SATs) for boundary and initial conditions provides a semi-discretization recipe for well posed linear initial boundary value problems (IBVPs). The SBP-SAT frame work in space guarantees semi-discrete energy stability combined with high order accuracy. 
Recent work have shown how to extend the linear semi-discrete SBP-SAT theory to nonlinear IBVPs. 
In this note we further 
extend the nonlinear semi-discrete energy stability to nonlinear fully discrete unconditional energy stability.  We prove that strong nonlinear semi-discrete energy stability, the orthogonality of time and space and a high order SBP-SAT discretization in time lead to an energy bound for the fully discrete approximation in terms of only externally given initial and boundary data.

\end{abstract}



\end{frontmatter}


\section{Introduction}
The SBP-SAT frame work \cite{svard2014review,fernandez2014review} guarantees linear semi-discrete energy stability for IBVPs and can be applied to most (if not all)
high order numerical discretization methods. 
Recent work \cite{nordstrom2022linear-nonlinear,NORDSTROM2024_BC,Nordström2025_open,Nordström2026_LIN,Nonlin_SBP}  have extended the 
SBP-SAT theory into the nonlinear regime 
which leads to provably stable semi-discrete approximations of nonlinear IBVPs as systems of nonlinear ordinary differential equations (ODEs).

ODEs 
have been dealt with previously. 
In strong stability preserving (SSP) methods \cite{SSP1,SSP3}, the higher order time-integration schemes are related to the explicit forward Euler scheme and aims at retaining its stability properties. B-stability \cite{B-stab1,B-stab3} is the nonlinear version of the linear A-stability concept and consider nonlinear test equations. Both these concepts are general in nature and not necessarily related to IBVPs.

Here we work directly on the semi-discretized IBVP and extend the nonlinear semi-discrete stability to nonlinear fully discrete unconditional stability. We prove
that strong nonlinear semi-discrete energy stability, the orthogonality of time and space and high order
SBP-SAT discretizations in time \cite{SBP_in_time_1,SBP_in_time_2,nordstrom2019} lead to an energy bound for the fully discrete IBVP approximation in terms of initial and boundary data only.



\section{Analysis}\label{sec:theory}
Consider the following general semi-discretized  nonlinear IBVP with a SAT term for boundary conditions
\begin{equation}\label{eq:nonlin_Semi}
P_V U_t + (Q(U)+R(U))U = SAT_S(U,G),  \quad t \geq 0
\end{equation}
augmented with the boundary conditions $B(U)-G=0$ imposed by the $SAT_S$ and an 
the initial condition $U-F=0$. 
In (\ref{eq:nonlin_Semi}),  $U(t)$ is the $n$ component long time-dependent 
solution discretized in N spatial nodes, $Q(U)$ is an almost skew-symmetric solution dependent matrix operator \cite{Nonlin_SBP}, 
the subscript $t$ indicate differentiation in time and the diagonal volume (or mass) matrix $P_V$ 
defines an $L_2$ similar energy norm (or semi-norm) $\|U\|^2_{P_V} = U^T P_V U$. For a dissipative IBVP, $R(U)$ is a positive dissipative symmetric 
matrix operator, otherwise it is zero. Furthermore,  the weak 
boundary conditions are implemented by the space penalty term $SAT_S(U,G)=\Sigma_S(U)(B(U)-G)$ where 
$\Sigma_S(U)$ is the space penalty matrix, $B(U)$ is the nonlinear boundary operator and $G$ the externally given 
boundary data.
\begin{remark}
The structure of the nonlinear system of ordinary differential equations (\ref{eq:nonlin_Semi}) in time was derived in \cite{nordstrom2022linear-nonlinear,NORDSTROM2024_BC,Nordström2025_open,Nordström2026_LIN,Nonlin_SBP}. It originates for example as a semi-discretization of nonlinear IBVPs related to the shallow water \cite{nordstrom2022linear-nonlinear}, Euler \cite{Nordstrom2022_Skew_Euler}, Navier-Stokes \cite{nordstrom2024skewsymmetric_jcp} or multi-phase flow equations \cite{Nordström2027}. 
\end{remark}

\subsection{Preliminaries}\label{Semi_discrete}
Next, we define semi-discrete energy stability (the key building block for fully discrete stability).
\begin{definition}
The semi-discretized IBVP (\ref{eq:nonlin_Semi}) is semi-discretely strongly energy stable if at time $t=T$, 
\begin{equation}\label{eq:semi-disc_stab}
\|U(T)\|^2_{P_V} \leq \|F\|^2_{P_V} + 2 \int\limits_{0}^T \|G\|^2_{P_B}dt
\end{equation}
holds. If (\ref{eq:semi-disc_stab}) holds for $G=0$ only, it is semi-discretely energy stable. In (\ref{eq:semi-disc_stab}) we introduced the $L_2$ similar boundary norm $\|G\|^2_{P_B}=G^T P_B G$ where the matrix $P_B$ is diagonal.
\end{definition}
We can now prove the following Lemma.
\begin{lemma}\label{lemma:semi-discrete}
Consider the semi-discretized IBVP (\ref{eq:nonlin_Semi}) and let $S(U,G)=(Q(U)+R)U-SAT_S(U,G)$ denote the complete spatial discretization vector. The IBVP (\ref{eq:nonlin_Semi}) is semi-discretely strongly energy stable if 
\begin{equation}\label{Bcon_est}
U^T S(U,G) \geq - G^T P_B G = - \|G\|^2_{P_B}
\end{equation}
holds. It is semi-discretely energy stable if (\ref{Bcon_est}) holds only for $G=0$. 
\end{lemma}
\begin{proof}
The discrete energy method (multiply 
(\ref{eq:nonlin_Semi}) with $2 U^T$ from the left) yields $\frac{d}{dt}\|U\|^2_{P_V} + 2 U^T S(U,G)=0$. Condition (\ref{Bcon_est}) and time integration leads directly to the semi-discrete estimate (\ref{eq:semi-disc_stab}).
\end{proof}
\begin{remark}\label{semi-fully-required} 
The semi-discrete estimate (\ref{Bcon_est}) holds for all time and will later be used at discrete times.
\end{remark}
\begin{remark}\label{semi-fully-required} 
The specific semi-discrete formulation  (\ref{eq:nonlin_Semi}) with  weak boundary conditions is not required. 
Any formulation $S(U,G)$ that produce the estimate (\ref{Bcon_est}) will be shown below to lead to fully discrete stability.
\end{remark}

\subsection{The fully discrete scheme using SBP-SAT in time}\label{fully_discrete}
By discretizing the system of nonlinear ODEs (\ref{eq:nonlin_Semi}) using SBP-SAT in time  we obtain 
\begin{equation}\label{eq:fully_disc_eq}
(D_T \otimes I_V) \vec U + (I_T \otimes P_V^{-1}) \vec S  = SAT_T(\vec U,\vec F)
\end{equation}
where $\vec U = (U_0, U_1,...,U_M)^T$ and $U_i$ denote $U(t_i)$ and similary with $\vec S$. To separate the dimensions space and time we have introduced the Kronecker symbol $\otimes$ \cite{horn2012}. The SBP difference operator in time is $D_T=P_T^{-1} Q_T$ and the  
time penalty term is $SAT_T(\vec U,\vec F)=(P_T^{-1} E_0 \otimes I_V)\Sigma_T(\vec U-\vec F)$.  
The 
SBP property \cite{svard2014review,fernandez2014review} 
is $Q_T + Q_T^T = B_T =  E_M - E_0$
where $E_0 = \rm{diag}(1,0,\dots,0)$ and $E_N = \rm{diag}(0,\dots,0,1)$. The diagonal matrix $P_T$ is positive definite and defines a quadrature formula. 
In $SAT_T$, $\Sigma_T$ is the time penalty matrix and $\vec F$ contain the initial data $F$ in the first component $F_0$ and arbitrary  values 
in the remaining components.

Next, we define fully discrete unconditional stability.
\begin{definition}
The discrete IBVP (\ref{eq:fully_disc_eq}) is discretely strongly unconditionally stable if the estimate 
\begin{equation}\label{eq:full-disc_stab_def}
\|U_M\|^2_{P_V} \leq \|F\|^2_{P_V} + 2\|\vec G\|^2_{P_T \otimes P_B}
\end{equation}
holds. If (\ref{eq:full-disc_stab_def}) holds for $\vec G=0$ only, it is discretely unconditionally stable. In (\ref{eq:full-disc_stab_def}) we introduced the $L_2$ similar time-boundary norm $\|\vec G\|^2_{P_T \otimes P_B}=\vec G^T (P_T \otimes P_B) \vec G$ where the matrix $(P_T \otimes P_B)$ is diagonal.
\end{definition}
\begin{remark} We will prove stability for one of the stages in the multistage version \cite{SBP_in_time_2} of SBP in time.
\end{remark}

The following Proposition is the main result of this note.
\begin{proposition}\label{prop_fully_discrete}
The discrete IBVP (\ref{eq:fully_disc_eq}) is unconditionally strongly  stable if (\ref{Bcon_est}) 
holds and $\Sigma_T=-I_V$.
\end{proposition}
\begin{proof}
The discrete energy method (multiply (\ref{eq:fully_disc_eq} with $2 \vec U^T(P_T \otimes P_V)$ from the left) yields 
\begin{equation}\label{eq:fully_disc_eq_proof1}
2\vec U^T (Q_T \otimes P_V) \vec U + 2\vec U^T (P_T \otimes I_V) \vec S  = 2\vec U^T (E_0 \otimes P_V) \Sigma_T(\vec U-\vec F).
\end{equation}
The first term on the lefthand side (LHS) in (\ref{eq:fully_disc_eq_proof1}) goes to
$2\vec U^T (Q_T \otimes P_V) \vec U = \vec U^T (Q_T + 
Q_T^T \otimes P_V) \vec U=$ 
\begin{equation}
\label{illustration1}
=\vec U^T (B_T \otimes P_V) \vec U =
{\begin{bmatrix}
      U_0 \\
      U_1 \\
    \vdots \\
      U_{M-1} \\
      U_M
\end{bmatrix}}^T
\begin{bmatrix}
       -P_V   &         &          &    &   \\
              &  0      &          &    &   \\
              &         &  \ddots  &    &   \\
              &         &          & 0  &   \\
              &         &          &    &  P_V
\end{bmatrix}
\begin{bmatrix}
      U_0 \\
      U_1 \\
    \vdots \\
     U_{M-1} \\
      U_M
\end{bmatrix}
=
U_M^T P_V U_M - U_0^T P_V U_0.
\end{equation}

The second term on the LHS in (\ref{eq:fully_disc_eq_proof1}) together with the notation $P_T=diag(P_0,P_1,...,P_{M-1},P_M)$ becomes
\begin{equation}
\label{illustration2}
2\vec U^T (P_T \otimes I_V) \vec S = 
{\begin{bmatrix}
      U_0 \\
      U_1 \\
    \vdots \\
      U_{M-1} \\
      U_M
\end{bmatrix}}^T
\begin{bmatrix}
       P_0 I_V  &         &          &    &   \\
                   & P_1 I_V      &          &    &   \\
              &         &  \ddots  &    &   \\
              &         &          & P_{M-1} I_V  &   \\
              &         &          &    &  P_M I_V
\end{bmatrix}
\begin{bmatrix}
      S_0 \\
      S_1 \\
    \vdots \\
     S_{M-1} \\
      S_M
\end{bmatrix}
= 2\sum_{i=0}^{i=M} P_i U_i^T  S_i.
\end{equation}
Next we recall that (\ref{Bcon_est}) in Lemma \ref{lemma:semi-discrete} implies that $U(t)^T S(t) \geq - G(t)^T P_B G(t)$ for all times and hence also for the discrete time levels in (\ref{illustration2}). This results in the estimate
\begin{equation}
\label{illustration3}
2\vec U^T (P_T \otimes I_V) \vec S = 2\sum_{i=0}^{i=M} P_i U_i^T  S_i \geq - 2\vec G^T (P_T \otimes P_B) \vec G = - 2\|\vec G\|^2_{P_T \otimes P_B}
\end{equation}
where $\vec G=(G_0,G_1,...,G_{M-1},G_M)^T$ includes all discrete time levels. 

The expanded SAT term on the RHS together with the choice $\Sigma_T=-I_V$ becomes
\begin{equation}
\label{illustration4}
2\vec U^T (E_0 \otimes P_V) \Sigma_T(\vec U-\vec F)= -2 U_0^T P_V U_0 + 2 U_0^T P_V F_0= -2 U_0^T P_V U_0 + 2 U_0^T P_V F.
\end{equation}
By adding and subtracting $F^T P_V F = \|F\|^2_{P_V}$ and using (\ref{illustration1})-(\ref{illustration4}) we obtain the final estimate
\begin{equation}\label{eq:full-disc_stab}
\|U_M\|^2_{P_V} \leq \|F\|^2_{P_V} + 2\|\vec G\|^2_{P_T \otimes P_B}-(U_0-F)^T P_V (U_0-F) \leq \|F\|^2_{P_V} + 2\|\vec G\|^2_{P_T \otimes P_B}
\end{equation}
which proves discrete unconditional strong stability.
\end{proof}
We also have the following Corollary which follows directly from the proof above.
\begin{corollary}\label{prop_fully_discrete}
The discrete IBVP (\ref{eq:fully_disc_eq}) is unconditionally  stable if (\ref{Bcon_est}) 
holds with $G=0$.
\end{corollary}
\begin{remark}\label{block-structure-SBP-in-time} 
The orthogonality of space and time enabled the tensor product form of the approximation which together with SBP-SAT in time isolated the semi-discrete stability condition in (\ref{Bcon_est}) and (\ref{illustration2}). 
\end{remark}

\section{Summary}\label{sec:conclusion}



We have shown how to extend nonlinear semi-discrete energy stable schemes for IBVPs to fully discrete energy stable schemes. 
We have proved that nonlinear semi-discrete strong energy stability, the orthogonality of time and space and high order accurate SBP-SAT discretizations in time lead to fully discrete nonlinear unconditional strong stability and energy bounds in terms of 
initial and boundary data only. 

\section*{Acknowledgment}
J.N. was supported by University of Johannesburg Global
Excellence and Stature Initiative Funding.


\bibliographystyle{elsarticle-num}
\bibliography{References_Jan}

\end{document}